\documentclass{amsart}

\usepackage{amssymb}
\usepackage{microtype}
\usepackage{fullpage}
\usepackage{hyperref}

\hypersetup{colorlinks=true, citecolor=blue, linkcolor=blue, urlcolor=blue, pdfstartview=FitH, pdfauthor=Antal Balog and Andras Biro and Gergely
Harcos and Peter Maga, pdftitle=The prime geodesic theorem in square mean}

\newcommand{\ZZ}{\mathbb{Z}}

\newcommand{\RR}{\mathbb{R}}
\newcommand{\CC}{\mathbb{C}}
\newcommand{\cH}{\mathcal{H}}

\newcommand{\eps}{\varepsilon}
\newcommand{\bs}{\backslash}
\renewcommand{\geq}{\geqslant}
\renewcommand{\leq}{\leqslant}

\newcommand{\PSL}{\mathrm{PSL}}

\DeclareMathOperator{\sym}{sym}

\theoremstyle{plain}

\newtheorem{theorem}{Theorem}

\theoremstyle{remark}
\newtheorem{remark}{Remark}

\theoremstyle{definition}

\begin{document}

\author{Antal Balog}
\author{Andr\'as Bir\'o}
\author{Gergely Harcos}
\author{P\'eter Maga}

\address{Alfr\'ed R\'enyi Institute of Mathematics, Hungarian Academy of Sciences, POB 127, Budapest H-1364, Hungary}\email{balog@renyi.hu, biroand@renyi.hu,gharcos@renyi.hu, magapeter@gmail.com}
\address{MTA R\'enyi Int\'ezet Lend\"ulet Automorphic Research Group}\email{balog@renyi.hu, biroand@renyi.hu, gharcos@renyi.hu, magapeter@gmail.com}
\address{Central European University, Nador u. 9, Budapest H-1051, Hungary}\email{harcosg@ceu.edu}

\title{The prime geodesic theorem in square mean}

\begin{abstract}
We strengthen the recent result of Cherubini--Guerreiro on the square mean of the error term in the prime geodesic theorem for $\PSL_2(\ZZ)$.
We also develop a short interval version of this result.
\end{abstract}

\subjclass[2010]{Primary 11F72; Secondary 11M36}

\keywords{Prime geodesic theorem, spectral exponential sums}

\thanks{Supported by NKFIH (National Research, Development and Innovation Office) grant K~119528 and by the MTA R\'enyi Int\'ezet Lend\"ulet Automorphic Research Group. Fourth author also supported by the Premium Postdoctoral Fellowship of the Hungarian Academy of Sciences.}

\maketitle

\section{Introduction}

The aim of this note is to provide a new upper bound for the square mean error in the classical prime geodesic theorem. For a brief introduction, let
$\Gamma:=\PSL_2(\ZZ)$ be the modular group, and let
\[\Psi_\Gamma(X):=\sum_{NP\leq X}\Lambda(P)\]
denote the usual Chebyshev-like counting function for the closed geodesics on the modular surface $\Gamma\bs\cH$. That is, $\log NP$ is the length of the closed geodesic $P$, and $\Lambda(P)=\log NP_0$ is the length of the underlying prime closed geodesic $P_0$. The closed geodesic $P$ (resp. $P_0$) is understood without orientation, hence it corresponds bijectively to an unordered pair of hyperbolic (resp. primitive hyperbolic) conjugacy classes in $\Gamma$ which are reciprocals of each other (cf.\ \cite{Sa1,Sa2}). In an original breakthrough, Iwaniec~\cite{Iw} proved that
\[\Psi_\Gamma(X)=X+O_\eps(X^{35/48+\eps})\]
for any $\eps>0$, the important point being that $35/48$ in the exponent is less than $3/4$. This constant was subsequently lowered to $7/10$ by Luo--Sarnak~\cite{LuSa}, $71/102$ by Cai~\cite{Ca}, and $25/36$ by Soundararajan--Young~\cite{SoYo}. For the last mentioned result, Balkanova--Frolenkov~\cite{BaFr} provided a new proof very recently. It is conjectured that the exponent $2/3+\eps$ or perhaps even $1/2+\eps$ is admissible (in which case it would be optimal). Our main result states that the exponent $7/12+\eps$ is valid in a square mean sense.
\begin{theorem}\label{thm1} Let $A>2$. Then, for any $\eps>0$ we have
\[\frac{1}{A}\int_A^{2A}\left|\Psi_\Gamma(X)-X\right|^2dX\ll_\eps A^{7/6+\eps}.\]
\end{theorem}
This estimate improves on the result of Cherubini--Guerreiro~\cite[Th.~1.4]{ChGu}, where the right hand side was $A^{5/4+\eps}$, and in fact our analysis is based on theirs. Incidentally, the exponents $7/12+\eps$ and $5/8+\eps$ also occur in the recent works of Petridis--Risager~\cite{PeRi} and Bir\'o~\cite{Bi} on the hyperbolic circle problem, although their averages are not fully analogous to ours.

Theorem~\ref{thm1} has the following simple consequence for short intervals. For $0\leq\eta\leq1$ we have
\[\frac1A\int_A^{2A}\left|\Psi_\Gamma(X)-\Psi_\Gamma(X-\eta X)-\eta X\right|^2dX\ll_\eps A^{7/6+\eps},\]
that is, the approximation $\Psi_\Gamma(X)-\Psi_\Gamma(X-\eta X)\approx\eta X$ is valid with error term $X^{7/12+\eps}$ in a square mean sense. For $\eta\geq A^{-1/6}$, this is the best we can say at the moment. However, for smaller $\eta$, we can obtain an improvement by tailoring our analysis to the present problem, with an average error term tending to $X^{1/2+\eps}$ as $\eta$ gets close to $A^{-1/2}$.
\begin{theorem}\label{thm3} Let $A>2$. Then, for any $\eps>0$ and $A^{-1/2}\log^2A\leq\eta<A^{-1/6}$ we have
\[\frac{1}{A}\int_A^{2A}\left|\Psi_\Gamma(X)-\Psi_\Gamma(X-\eta X)-\eta X\right|^2dX\ll_\eps A^{5/4+\eps}\eta^{1/2}.\]
\end{theorem}
\begin{remark} Theorem~\ref{thm3} can be improved for very small $\eta$ by employing \cite[Th.~8.3]{BaFr}.
Specifically, on the right hand side of the bound, $A^{1+\eps}$ is admissible for $A^{-1/2}\log^2A\leq\eta<A^{-4/9}$, and $A^{5/3+\eps}\eta^{3/2}$ is admissible for $A^{-4/9}\leq\eta<A^{-5/12}$. See also Remark~\ref{remark2} below Theorem~\ref{thm2}.
\end{remark}

The paper is structured as follows. The overall strategy is already present in Iwaniec's seminal paper \cite{Iw}, but we also rely crucially on the work of Cherubini--Guerreiro~\cite{ChGu} and Luo--Sarnak~\cite{LuSa}. In Section~\ref{section4}, we reduce Theorems~\ref{thm1} and \ref{thm3} to the estimation of a certain spectral exponential sum. This reduction ultimately follows from Selberg's trace formula, although we do not invoke it explicitly. In Section~\ref{section3}, we prove Theorem~\ref{thm2}, which contains the necessary bounds for the spectral exponential sum. This proof is ultimately an application of Kuznetsov's trace formula, which again remains in the background, however. Section~\ref{section2} prepares the scene, incorporating a key idea of Iwaniec~\cite{Iw}.

\section{Reduction to Kuznetsov's trace formula}\label{section2}

Let $\{u_j\}$ be an orthonormal Hecke eigenbasis of the space of Maass cusp forms on $\Gamma\bs\cH$. Denoting by $1/4+t_j^2$ the Laplace eigenvalue of $u_j$ with the sign convention $t_j>0$, we have the Fourier decomposition
\[u_j(x+iy)=\sqrt{y}\sum_{n\neq 0}\rho_j(n)K_{it_j}(2\pi|n|y)e(nx).\]
The Fourier coefficients $\rho_j(n)$ are proportional to the Hecke eigenvalues $\lambda_j(n)$,
\begin{equation}\label{proportional}
\rho_j(n)=\rho_j(1)\lambda_j(n).
\end{equation}
The Hecke eigenvalues are real, and they satisfy the multiplicativity relations
\[\lambda_j(m)\lambda_j(n)=\sum_{d\mid\gcd(m,n)}\lambda_j\left(\frac{mn}{d^2}\right).\]
In particular, the symmetric square $L$-function of $u_j$ satisfies
\begin{equation}\label{symsquareLfunction}
L(s,\sym^2 u_j)=\zeta(2s)\sum_{n=1}^\infty\frac{\lambda_j(n^2)}{n^s}=\frac{\zeta(2s)}{\zeta(s)}\sum_{n=1}^\infty\frac{\lambda_j(n)^2}{n^s},
\qquad\Re s>1,\end{equation}
in the region of absolute convergence of both Dirichlet series.

Concerning the distribution of Laplace eigenvalues, we record Weyl's law as (cf.\ \cite[(11.5)]{Iw2})
\begin{equation}\label{Weyllaw}
\#\{j: t_j\leq T\}=\frac{T^2}{12}+O(T\log T).
\end{equation}
In fact a finer asymptotic expansion is available, see \cite[Ch.~11, (2.12)]{He} or \cite[Th.~7.3]{Ve}.

With Kuznetsov's trace formula in mind, we introduce the harmonic weights
\begin{equation}\label{harmonicweights}
\alpha_j:=\frac{|\rho_j(1)|^2}{\cosh(\pi t_j)}=\frac{2}{L(1,\sym^2 u_j)},
\end{equation}
which by \cite[Th.~8.3]{Iw2} and \cite[Th.~0.2]{HoLo} satisfy the convenient bounds
\begin{equation}\label{harmonicweightbounds}
t_j^{-\eps}\ll_\eps\alpha_j\ll_\eps t_j^\eps.
\end{equation}
For arbitrary $X,T>2$, we borrow from \cite[Lemma~7]{DeIw} the test function (see also \cite[p.~234]{LuSa} and \cite[Lemma~2.2]{BaFr})
\[\varphi(x):=\frac{\sinh\beta}{\pi}x\exp(ix\cosh\beta)\qquad\text{with}\qquad \beta:=\frac{\log X}{2}+\frac{i}{2T},\]
whose Bessel transform
\[\hat\varphi(t):=\frac{\pi i}{2\sinh(\pi t)}\int_0^\infty\bigl(J_{2it}(x)-J_{-2it}(x)\bigr)\varphi(x)\frac{dx}{x}\]
satisfies
\begin{equation}\label{Besseltransform}
\hat\varphi(t)=\frac{\sinh(\pi t+2\beta it)}{\sinh(\pi t)}=X^{it}e^{-t/T}+O\left(e^{-\pi t}\right).
\end{equation}

Following \cite{Iw,LuSa}, we consider the spectral-arithmetic average (cf.\ \eqref{symsquareLfunction})
\[\sum_j\alpha_j\hat\varphi(t_j)\sum_n h(n)\lambda_j(n)^2=
  \sum_j\alpha_j\hat\varphi(t_j)\frac{1}{2\pi i}\int_{(2)}\tilde h(s)\frac{\zeta(s)}{\zeta(2s)}L(s,\sym^2 u_j)\,ds,\]
where $h:(0,\infty)\to\RR$ is a smooth compactly supported function with holomorphic Mellin transform $\tilde h:\CC\to\CC$.
We choose $h$ such that it is supported in some dyadic interval $[N,2N]$ for $N>1$, and it satisfies $h^{(j)}\ll_jN^{-j}$ and $\tilde h(1)=N$. Then also
\begin{equation}\label{Mellinbound}
\tilde h(s)=\frac{(-1)^j}{s(s+1)\dots (s+j-1)}\int_0^\infty h^{(j)}(x)\, x^{s+j}\,\frac{dx}{x}\ll_{\sigma,j}\frac{N^{\sigma}}{(1+|s|)^j},\qquad\Re(s)=\sigma,
\end{equation}
where the implied constant depends continuously on $\sigma$. More precisely, the identity is meant for $s$ outside $\{0,-1,-2,\dots\}$, but the inequality holds even at these exceptional points. Shifting the contour, we obtain by the residue theorem and \eqref{harmonicweights},
\[\sum_j\alpha_j\hat\varphi(t_j)\sum_n h(n)\lambda_j(n)^2=\frac{12N}{\pi^2}\sum_j\hat\varphi(t_j)+
  \sum_j\alpha_j\hat\varphi(t_j)\frac{1}{2\pi i}\int_{(1/2)}\tilde h(s)\frac{\zeta(s)}{\zeta(2s)}L(s,\sym^2 u_j)\,ds.\]
Using also the approximation \eqref{Besseltransform}, we obtain after some rearrangement,
\begin{equation}\begin{split}\label{startingpoint}
\sum_j X^{it_j}e^{-t_j/T}=&\ O(1)+\frac{\pi^2}{12N}\sum_n h(n)\sum_j\alpha_j\hat\varphi(t_j)\lambda_j(n)^2\\
&\ -\frac{\pi^2}{12N}\,\frac{1}{2\pi i}\int_{(1/2)}\tilde h(s)\frac{\zeta(s)}{\zeta(2s)}\sum_j\alpha_j\hat\varphi(t_j)L(s,\sym^2 u_j)\,ds.
\end{split}\end{equation}
This formula is equivalent to \cite[(3.8)]{BaFr}, and we have included the proof for the sake of completeness.
We stress that the spectral weights $\hat\varphi(t_j)$ depend on the parameters $X,T>2$.

\section{Spectral exponential sums in square mean}\label{section3}

We shall estimate the spectral exponential sum \eqref{startingpoint}, in square mean over $A\leq X\leq 2A$,
by combining \eqref{startingpoint} with the analysis of Cherubini--Guerreiro~\cite{ChGu} and Luo--Sarnak~\cite{LuSa}.
Specifically, on the right hand side of \eqref{startingpoint}, the square mean of the first $j$-sum can be estimated via Kuznetsov's formula and the Hardy--Littlewood--P\'olya inequality (cf.\ \cite[Lemma~4.2]{ChGu}), while the square mean of the second $j$-sum can be estimated in terms of the spectral second moment of symmetric square $L$-functions (cf.\ \cite[(33)]{LuSa}). This way we obtain the following improvement of \cite[Prop.~4.5]{ChGu}.
\begin{theorem}\label{thm2} Let $A>2$. Then, for any $\eps>0$ we have
\begin{equation}\label{bound2}\frac{1}{A}\int_A^{2A}\Biggl|\sum_{t_j\leq T} X^{it_j}\Biggr|^2dX
\ll_\eps(AT)^\eps\begin{cases}
T^3,& 0<T\leq A^{1/6};\\
A^{1/4}T^{3/2},& A^{1/6}<T\leq A^{1/2};\\
T^2,&A^{1/2}<T.\end{cases}
\end{equation}
In particular, the left hand side can always be bounded as $\ll_\eps(AT)^\eps A^{1/6}T^2$.
\end{theorem}
\begin{remark}\label{remark2} Theorem~\ref{thm2} can be refined in the medium range by employing \cite[Th.~8.3]{BaFr}. Specifically, $A^{1/4}T^{3/2}$ can be improved to $A^{1/2+\theta}T^{1/2}$ for $A^{1/4+\theta}<T\leq A^{1/3+2\theta/3}$, and to $T^2$ for $A^{1/3+2\theta/3}<T\leq A^{1/2}$. Note that for $\theta$ any value exceeding $1/6$ is admissible by the celebrated work of Conrey--Iwaniec~\cite[Cor.~1.5]{CoIw}.
\end{remark}

Following the proof of \cite[Prop.~4.5]{ChGu}, which is based on \cite[pp.~235--236]{LuSa}, we see that \eqref{bound2} can be deduced from the following smoothened variant, itself a strengthening of \cite[Lemma~4.4]{ChGu}:
\begin{equation}\label{bound4}\frac{1}{A}\int_A^{2A}\Biggl|\sum_j X^{it_j}e^{-t_j/T}\Biggr|^2dX
\ll_\eps(AT)^\eps\begin{cases}
T^3,& 0<T\leq A^{1/6};\\
A^{1/4}T^{3/2},& A^{1/6}<T\leq A^{1/2};\\
T^2,&A^{1/2}<T.\end{cases}
\end{equation}
We shall assume here that $T>2$, since otherwise \eqref{bound4} is trivial. As a first step for the proof of \eqref{bound4}, we change in \eqref{startingpoint} the second occurrence of $\hat\varphi(t_j)$ to $X^{it_j}e^{-t_j/T}$, and we restrict the integration to $|\Im(s)|\leq T^\eps$.
The error resulting from this change is $O_\eps(1)$ by \eqref{Mellinbound} and standard bounds for the symmetric square $L$-function and the Riemann zeta function. Then, applying the Cauchy--Schwarz inequality multiple times and standard bounds for the Riemann zeta function, we arrive at
\begin{align*}&\Biggl|\sum_j X^{it_j}e^{-t_j/T}\Biggr|^2\ll_\eps 1+
\frac{1}{N}\sum_{N\leq n\leq 2N}\Biggl|\sum_j\alpha_j\hat\varphi(t_j)\lambda_j(n)^2\Biggr|^2\\
&+\frac{T^\eps}{N}\int_{-T^\eps}^{T^\eps}\Biggl|\sum_j\alpha_jX^{it_j}e^{-t_j/T}L(1/2+i\tau,\sym^2 u_j)\Biggr|^2d\tau.
\end{align*}
We abbreviate
\[L_j(\tau):=L(1/2+i\tau,\sym^2 u_j),\]
and we average over $A\leq X\leq 2A$. Applying \cite[Lemma~4.2]{ChGu}\footnote{The definitions of $\nu_j(n)$ and $\hat\phi(t)$ in \cite{ChGu,LuSa} are slightly in error, in particular their $\rho_j(n)=\nu_j(n)\cosh(\pi t_j/2)$ should really be
$\rho_j(n)=\nu_j(n)\cosh(\pi t_j)^{1/2}$. With this correction, $|\nu_j(n)|^2$ in \cite{ChGu,LuSa} agrees with our $\alpha_j\lambda_j(n)^2$, thanks to \eqref{proportional} and \eqref{harmonicweights}. For precise versions of the relevant Kuznetsov formula, see \cite[Th.~2]{Ku} and \cite[Th.~9.5]{Iw2}.} for the contribution of the $n$-sum on the right hand side, we obtain
\begin{align*}&\frac{1}{A}\int_A^{2A}\Biggl|\sum_j X^{it_j}e^{-t_j/T}\Biggr|^2dX\ll_\eps (NA^{1/2}+T^2)(ANT)^\eps\\
&+\frac{T^\eps}{N}\int_{-T^\eps}^{T^\eps}\frac{1}{A}\int_A^{2A}\Biggl|\sum_j\alpha_jX^{it_j}e^{-t_j/T}L_j(\tau)\Biggr|^2 dX\,d\tau.
\end{align*}

We apply the Cauchy--Schwarz inequality one more time to facilitate the upcoming analysis. Specifically, we distribute the spectral parameters $t_j$ on the right hand side into intervals of length $T$, and this way we get
\[\Biggl|\sum_j\alpha_jX^{it_j}e^{-t_j/T}L_j(\tau)\Biggr|^2
\ll\sum_{m=1}^\infty m^2\,\Biggl|\sum_{(m-1)T\leq t_j<mT}\alpha_jX^{it_j}e^{-t_j/T}L_j(\tau)\Biggr|^2.\]
Therefore, with the notation
\begin{equation}\label{Idef}
I(T,A,m,\tau):=\frac{1}{A}\int_A^{2A}\Biggl|\sum_{(m-1)T\leq t_j<mT}\alpha_jX^{it_j}e^{-t_j/T}L_j(\tau)\Biggr|^2dX,
\end{equation}
we infer
\begin{equation}\label{bound1}
\frac{1}{A}\int_A^{2A}\Biggl|\sum_j X^{it_j}e^{-t_j/T}\Biggr|^2dX\ll_\eps
(NA^{1/2}+T^2)(ANT)^\eps+\frac{T^\eps}{N}\sup_{|\tau|\leq T^\eps}\sum_{m=1}^\infty m^2\,I(T,A,m,\tau).
\end{equation}
We bound $I(T,A,m,\tau)$ by squaring out the $j$-sum in \eqref{Idef}, then integrating explicitly in $X$, and finally using \eqref{harmonicweightbounds} for $\alpha_j$:
\begin{align*}I(T,A,m,\tau)&\ll_\eps T^\eps e^{-2m}\sum_{(m-1)T\leq t_j,t_k<mT}\frac{\left|L_j(\tau)L_k(\tau))\right|}{1+|t_j-t_k|}\\
&\leq\ \:\frac{T^\eps e^{-2m}}{2}\sum_{(m-1)T\leq t_j,t_k<mT}\frac{\left|L_j(\tau)\right|^2+\left|L_k(\tau)\right|^2}{1+|t_j-t_k|}\\
&=\ \ T^\eps e^{-2m}\sum_{(m-1)T\leq t_j<mT}\left|L_j(\tau)\right|^2\sum_{(m-1)T\leq t_k<mT}\frac{1}{1+|t_j-t_k|}.
\end{align*}
By the Weyl law \eqref{Weyllaw}, the last $k$-sum is
\begin{equation}\label{bound3}
\sum_{(m-1)T\leq t_k<mT}\frac{1}{1+|t_j-t_k|}\ \leq\ \sum_{\ell=1}^{\lceil T\rceil}\frac{1}{\ell}\sum_{\substack{(m-1)T\leq t_k<mT\\\ell-1\leq |t_j-t_k|<\ell}}1\ \ll_\eps\ (mT)^{1+\eps}\sum_{\ell=1}^{\lceil T\rceil}\frac{1}{\ell}\ \ll_\eps(mT)^{1+2\eps},
\end{equation}
whence
\[I(T,A,m,\tau)\ll_\eps (mT)^{1+\eps}\,e^{-2m}\sum_{(m-1)T\leq t_j<mT}\left|L_j(\tau)\right|^2.\]
For the last sum, we apply the spectral second moment bound of Luo--Sarnak \cite[(33)]{LuSa}, obtaining
\[I(T,A,m,\tau)\ll_\eps (mT)^{3+\eps}(1+|\tau|)^{5+\eps}\,e^{-2m}.\]
In combination with \eqref{bound1}, this yields
\[\frac{1}{A}\int_A^{2A}\Biggl|\sum_j X^{it_j}e^{-t_j/T}\Biggr|^2dX\ll_\eps(NA^{1/2}+T^2)(ANT)^\eps+\frac{T^{3+\eps}}{N}.\]

The last bound improves on the display before \cite[Prop.~4.5]{ChGu} in that we have $T^{3+\eps}$ in place of $T^{4+\eps}$.
We optimize by setting $N:=A^{-1/4}T^{3/2}$, which exceeds $1$ if and only if $T>A^{1/6}$. Assuming this, we obtain \eqref{bound4} readily. For $T\leq A^{1/6}$ we estimate the left hand side of \eqref{bound4} more directly but along the same ideas. Specifically, let us distribute the spectral parameters $t_j$ into intervals of length $T$ as before, apply the Cauchy--Schwarz inequality for the resulting $m$-sum, square out the various $j$-subsums, integrate explicitly in $X$, and then apply the Weyl law \eqref{Weyllaw}. We obtain (cf.\ \eqref{bound3})
\begin{align*}\frac{1}{A}\int_A^{2A}\Biggl|\sum_j X^{it_j}e^{-t_j/T}\Biggr|^2dX
&\ll\:\:\sum_{m=1}^\infty m^2 e^{-2m}\sum_{(m-1)T\leq t_j,t_k<mT}\frac{1}{1+|t_j-t_k|}\\
&\ll_\eps\sum_{m=1}^\infty m^2 e^{-2m}(m^{1+\eps}T^2)(mT)^{1+\eps}\ll_\eps T^{3+\eps},
\end{align*}
which is equivalent to \eqref{bound4} for $T\leq A^{1/6}$. The proof of Theorem~\ref{thm2} is complete.

\section{Prime geodesic error terms in square mean}\label{section4}

In this section, we deduce Theorem~\ref{thm1} and \ref{thm3} from Theorem~\ref{thm2}. For both theorems, we shall assume (without loss of generality) that $A>100$.

Our deduction of Theorem~\ref{thm1} follows almost verbatim the argument of Cherubini--Guerreiro right after the proof of \cite[Prop.~4.5]{ChGu}. We reproduce this argument (with small corrections), because we shall use certain steps from it in the proof of Theorem~\ref{thm3}.  Our starting point is the explicit formula for $\Psi_\Gamma(X)$ established by Iwaniec~\cite[Lemma~1]{Iw},
\begin{equation}\label{explicit}
\Psi_\Gamma(X)=X+\sum_{|t_j|\leq T}\frac{X^{1/2+it_j}}{1/2+it_j}+O\left(\frac XT\log^2X\right),\qquad 2<T\leq\frac{X^{1/2}}{\log^2X}.
\end{equation}
Here the notation $|t_j|\leq T$ means that the sum runs through the spectral parameters $\pm t_j$ with $t_j\leq T$ (recall our sign convention $t_j>0$). With the notation
\[R(X,T):=\sum_{t_j\leq T} X^{it_j},\]
the spectral sum in the explicit formula can be expressed as twice the real part of
\[\sum_{t_j\leq T}\frac{X^{1/2+it_j}}{1/2+it_j}=X^{1/2}\frac{R(X,T)}{1/2+iT}+iX^{1/2}\int_1^T\frac{R(X,U)}{(1/2+iU)^2}\;dU.\]
We specify $T:=A^{1/2}/\log^2A$, and we note that $T\leq X^{1/2}/\log^2X$ holds for any $X\geq A$ by the assumption $A>100$.
Applying the Cauchy--Schwarz inequality several times,
\begin{align*}
\frac{1}{A}\int_A^{2A}\Biggl|\sum_{t_j\leq T}\frac{X^{1/2+it_j}}{1/2+it_j}\Biggr|^2dX
&\ll\int_A^{2A}\left|\frac{R(X,T)}{1/2+iT}\right|^2dX+\int_A^{2A}\left|\int_1^T\frac{R(X,U)}{(1/2+iU)^2}\;dU\right|^2dX\\
&\ll\frac{1}{T^2}\int_A^{2A}|R(X,T)|^2dX+\log T\int_1^T\left(\int_A^{2A}|R(X,U)|^2 dX\right)\frac{dU}{U^3}.
\end{align*}
On the right hand side, the first term is $O_\eps(A^{1+\eps})$ and the second term is $O_\eps(A^{7/6+\eps})$
by Theorem~\ref{thm2}. Noting also that the error term in \eqref{explicit} is $O_\eps(A^{1/2+\eps})$, we obtain the bound in Theorem~\ref{thm1}.

Now we prove Theorem~\ref{thm3}. We specify $T:=A^{1/2}/\log^2A$ as before. The condition $A^{-1/2}\log^2A\leq\eta<A^{-1/6}$ then yields $T^{-1}\leq\eta<1/2$. By the explicit formula \eqref{explicit},
\[\Psi_\Gamma(X)-\Psi_\Gamma(X-\eta X)-\eta X=\sum_{|t_j|\leq T}X^{1/2+it_j}\frac{1-(1-\eta)^{1/2+it_j}}{1/2+it_j}+O_\eps(A^{1/2+\eps}),\]
and we need to estimate the square mean of this expression over $A\leq X\leq 2A$. It suffices to do this with the restriction $t_j>0$ on the right hand side, since the original sum over $|t_j|\leq T$ is twice the real part of the new sum over $t_j\leq T$. The contribution of the spectral parameters $t_j\leq 1/\eta$ can be rewritten and bounded by the Cauchy--Schwarz inequality as
\[\frac{1}{A}\int_A^{2A}\Biggl|\int_{1-\eta}^1\left(\frac X\xi\right)^{1/2}R(X\xi,1/\eta)\,d\xi\Biggr|^2dX
\leq\left(\int_{1-\eta}^1\frac{d\xi}{\xi}\right)\left(\int_{1-\eta}^1\int_A^{2A}|R(X\xi,1/\eta)|^2\,dX\,d\xi\right).\]
The $X$-integral is $O_\eps(A^{5/4+\eps}\eta^{-3/2})$ by Theorem~\ref{thm2}, hence the right hand side is $O_\eps(A^{5/4+\eps}\eta^{1/2})$. The contribution of the spectral parameters $1/\eta<t_j\leq T$ is bounded by
\begin{equation}\label{twointegrals}
\frac{1}{A}\int_A^{2A}\Biggl|\sum_{1/\eta<t_j\leq T}\frac{X^{1/2+it_j}}{1/2+it_j}\Biggr|^2dX+
\frac{1}{B}\int_B^{2B}\Biggl|\sum_{1/\eta<t_j\leq T}\frac{X^{1/2+it_j}}{1/2+it_j}\Biggr|^2dX,
\end{equation}
where $B$ abbreviates $(1-\eta)A$. These integrals are very similar to the one we encountered in the proof of Theorem~\ref{thm1}, so we can be brief. The first integral in \eqref{twointegrals} can be bounded by partial summation, the Cauchy--Schwarz inequality, and Theorem~\ref{thm2} as
\begin{align*}
&\ll\:\:\eta^2\int_A^{2A}|R(X,1/\eta)|^2\,dX+\frac{1}{T^2}\int_A^{2A}|R(X,T)|^2\,dX+
\log(\eta T)\int_{1/\eta}^T\left(\int_A^{2A}|R(X,U)|^2\,dX\right)\frac{dU}{U^3}\\
&\ll_\eps A^{5/4+\eps}\eta^{1/2}+A^{1+\eps}+A^{5/4+\eps}\int_{1/\eta}^T U^{-3/2}\,dU\ll A^{5/4+\eps}\eta^{1/2}.
\end{align*}
Similarly, the second integral in \eqref{twointegrals} is $O_\eps(B^{5/4+\eps}\eta^{1/2})$, hence also $O_\eps(A^{5/4+\eps}\eta^{1/2})$. Finally, the contribution of the error term in \eqref{explicit} is $O_\eps(A^{1+\eps})$. The proof of Theorem~\ref{thm3} is complete.

\end{document}